\input vanilla.sty
\scaletype{\magstep1}
\scalelinespacing{\magstep1}
\def\bull{\vrule height .9ex width .8ex depth -.1ex}


\pageno=-1

\title Isometric stability property of certain Banach spaces
\endtitle

\author Alexander Koldobsky
\endauthor 

\vskip1truecm

\subheading{Abstract} Let $E$ be one of the spaces $C(K)$ and $L_1$,
$F$ be an arbitrary Banach space, $p>1,$ and  $(X,\sigma)$ be a 
space with a finite measure. We prove that $E$ is isometric to a
subspace of the Lebesgue-Bochner space $L_p(X;F)$ only if $E$ is
isometric to a subspace of $F.$ Moreover, every isometry
 $T$ from $E$ into $L_p(X;F)$ has the form
$Te(x)=h(x)U(x)e, e\in E,$  where $h:X\rightarrow R$
is a measurable function and, for every $x\in X,$ $U(x)$ is
an isometry from $E$ to $F.$   

\vskip1truecm

\subheading{Mailing address:}

Till August 1,1993: 

Department of Mathematics

University of Missouri-Columbia

Columbia, MO 65211, USA

\vskip1truecm

After August 1,1993:

Division of Mathematics, Computer Science and Statistics

University of Texas at San Antonio

San Antonio, TX 78249, USA 

\vskip1truecm

\subheading{AMS classification:} 46B04, 47B80

\subheading{Key words:} isometries, Lebesgue-Bochner spaces, 
random operators

\subheading{Running title:} Isometric embeddings 

\newpage
\pageno=1

\title Isometric stability property of certain Banach spaces
\endtitle

\author Alexander Koldobsky\footnote{After August 1, 1993 the author's
address will be:Division of Mathematics, Computer Science
and Statistics, University of Texas at San Antonio,
San Antonio, TX 78249.}
\\Department of Mathematics\\ University of Missouri-Columbia\\
Columbia, MO 65211\\ 
\endauthor 

\vskip1truecm

\subheading{1.Introduction}

Let $E$ and $F$ be Banach spaces, $p\geq 1,$ $(X,\sigma)$
be a finite measure space, and $L_p(X,F)$ be the Lebesque-Bochner space
of (equivalence classes of) strongly measurable functions $f:X\mapsto F$ 
with 
$$\|f\|^p = \int_X \|f(x)\|^p d\sigma(x) < \infty.$$

We show that if $E=C(K)$ (with $K$ being a compact metric space) 
or $E=L_1$ then the space $E$ can be 
isometric to a subspace of $L_p(X,F)$ with $p>1$ only if $E$ is 
isometric to a subspace of $F.$

The isomorphic version of this result has been proved
for the spaces $E=c_0$ (S.Kwapien [5] and J.Bourgain [1]),
$E=l_1$ (G.Pisier [7]) and $E=l_{\infty}$ (J.Mendoza [6], see also
[2] for generalizations to Kothe spaces of vector valued functions). 
E.Saab and P.Saab [9] have proved
the isomorphic version for the space $E=L_1$ under the
additional assumption that $F$ is a dual space.   

For any $1\leq p\leq q \leq r,$ the space $L_q$ is isometric to a subspace 
of $L_p(X,L_r)$ (see [8]). On the other hand, if $r>2, r\neq q, q\neq 2$ the
space $L_q$ is not isomorphic to a subspace of $L_r.$ Thus, 
both isometric and isomorphic versions fail to be true
in the case where $E=L_q, q>1.$ (The space $L_2$ is isometric to a subspace
of $L_p$ for any $p>0,$ so it is isometric to a subspace of $L_p(X,F)$
for any space $F)$

Besides proving the above mentioned result 
for the spaces $C(K)$ and $L_1,$ we completely characterize
isometric embeddings of these spaces into $L_p$-spaces of vector 
valued functions.

Denote by $I(E,F)$ the set of isometries from $E$ to $F.$
A mapping $U:X\mapsto I(E,F)$ is called strongly measurable
if, for each $e\in E,$ the function $\|U(x)e\|$ is measurable
on $X.$ If the set $I(E,F)$ is non-empty then, obviously,
for every strongly measurable mapping $U:X\mapsto I(E,F)$ and
every function $h:X\mapsto R$ with $\|h\|_{L_p(X)}=1,$
the operator $T:E\mapsto L_p(X,F)$ defined by
$$Te(x)=h(x)U(x)e, e\in E \tag{1}$$
is an isometry.

We prove below that, for $E=C(K)$ or $E=L_1$ and for an arbitrary
space $F,$ every isometry from $E$ to $L_p(X,F)$ has the form (1).
The question if all isometries from $E$ to $L_p(X,F)$ have the form
(1) makes sense and has some applications even if $E=F.$
For example, if $E=F=L_q, q>1, p>1,$ then every isometry from $L_q$ 
to $L_p(X,L_q)$  has the form (1) if and only if $p\neq q, q\neq 2, 
q\not\in (p,2).$ This result was proved and applied to the 
description of isometries of Lebesgue-Bochner spaces in the paper [3].

As a consequence of the characterization of isometries from 
$E=C(K)$ or $E=L_1$ to $L_p(X,F)$ we obtain the following result on random
operators. Suppose that $(X,\sigma)$ is a probability space, let
$p>1,$ $L(E,F)$ be the space of linear operators from $E$ to $F$
and consider a random operator $S:X\mapsto L(E,F)$ which is an isometry
in average, i.e. $S$ is strongly measurable and, for every $e\in E,$
$$\|e\|^p = \int_X \|S(x)e\|^p d\sigma(x).$$
Then the random operator $S$ is an isometry (up to a constant) with
probability 1 (namely, operators $S(x)$ are isometries multiplied
by constants for almost all $x\in X).$ In fact, if we define an
operator $T:E\mapsto L_p(X,F)$ by $Te(x)=S(x)e$ then $T$ is an isometry
and (1) implies the desired result. In the case $E=C(K),$
the result about random operators has been proved before [4] under
the assumption that the operators $S(x)$ are bounded.

\subheading{2.Main results}

The proofs are based on the following simple fact.

\proclaim{Lemma 1} Let $E, F$ be Banach spaces, $p>1,$ $(X,\sigma)$ 
be a finite measure space and $T$ be an isometry from $E$ to $L_p(X,F).$
If $e,f$ are elements from $E$ with $\|e\|=\|f\|=1$ and  having
a common tangent functional (i.e. there exists $x^{*}\in E^{*}$ with
$\|x^{*}\|=1, x^{*}(e)=x^{*}(f)=1)$ then, for 
almost all (with respect to $\sigma)$ $x\in X,$ we have 

(i) $\|Te(x)\|=\|Tf(x)\|,$ 

(ii) for every $\alpha>0,$ $\|Te(x)+\alpha Tf(x)\|=\|Te(x)\|+ 
\alpha \|Tf(x)\|.$  
\endproclaim

\demo{Proof} Obviously, $\|e+\alpha f\|=1+\alpha$ for every $\alpha>0,$ 
and we have

$$(1+\alpha)^p = \|e+\alpha f\|^{p}=\int_{X} \|Te(x)+\alpha Tf(x)\|^{p} 
d\sigma (x) \leq$$ 
$$\int_{X} (\|Te(x)\|+\alpha \|Tf(x)\|)^{p} d\sigma (x) \tag{2}$$ 

For $\alpha=0,$ (2) turns into an equality. Therefore, we get 
a correct inequality if we
take in both sides of (2) the right-hand derivatives
at the point $\alpha=0$ and apply Holder's inequality:

$$1 \leq \int_{X} \|Te(x)\|^{p-1} \|Tf(x)\| d\sigma (x) \leq$$

$$ \Bigl( \int_{X} \|Te(x)\|^{p} d\sigma (x)\Bigr)^{p-1\over p}
\Bigl( \int_{X} \|Tf(x)\|^{p} d\sigma (x)\Bigr)^{1/p} = 1 $$

>From the conditions for equality in Holder's inequality, we conclude
that, for almost all (with respect to $\sigma)$  $x\in X,$
$\|Te(x)\|= c\|Tf(x)\|$ where $c$ is a constant. Further,
$$1=\int_{X} \|Te(x)\|^{p} d\sigma (x)=
\int_{X} c^{p}\|Tf(x)\|^{p} d\sigma (x)=c^{p}$$
and, hence, $c=1.$

It is clear now that (2) is, in fact, an equality. Hence, for almost
all $x\in X,$ we have 

$$\|Te(x)+\alpha Tf(x)\|=\|Te(x)\|+ \alpha \|Tf(x)\| \tag{3}$$

for every $\alpha>0$ and the proof is complete. \bull \enddemo

Now we are able to prove the main result.

\proclaim{Theorem 1} Let $p>1,$ $K$ be a compact metric space,
$(X,\sigma), (Y,\nu)$ be spaces with finite measures, and $F$ be 
an arbitrary Banach space. Let $E$ be either the space 
$L_1=L_1(Y,\nu)$ or
any subspace of $C(K)$ containing the function $1(k)\equiv 1$. Then

(i) If  $E$ is isometric to a subspace of $L_p(X;F)$   then
$E$ is isometric to a subspace of $F$ and the set $I(E,F)$
is non-empty.

(ii) If $T$ is an isometry from $E$ into $L_p(X;F)$  
then there exist
a measurable function $h:X\rightarrow R$  and a strongly 
measurable mapping $U:X\rightarrow I(E,F)$ such that
 $Te(x)=h(x)U(x)e$ for every $e\in E.$ \endproclaim 

\demo{Proof} We start with the case $E=L_1.$

Any two functions $e$ and $f$ from $L_1$ 
with disjoint supports in $Y$ have a common tangent functional
so we can apply Lemma 1 to any pair of normalized functions 
with disjoint supports. 

Decompose the set $Y$ into two parts $Y=Y_1\cup Y_2, Y_1\cap Y_2
=\emptyset, \nu (Y_i)>0, i=1,2.$ Fix a function $e_0\in L_1(Y_1),
\|e_0\|=1$ and put $h(x)=\|Te_0(x)\|, x\in X.$

Let $f_k, k\in N$ be a sequence of linearly independent functions
with supports in $Y_2$ such that their linear span is dense in 
$L_1(Y_2).$ Denote 
by $D$ the set of linear combinations of functions $f_k$ with
rational coefficients. Given fixed representatives $Tf_k$ from
the corresponding equivalence classes of functions from 
the space $L_p(X;F),$ define an operator $T(x):D\mapsto F$ for every 
$x\in X$ by $T(x)(\sum \lambda_{i} f_i)= \sum \lambda_{i} 
Tf_{i}(x), \lambda_i \in Q.$

It follows from the statement (i) of Lemma 1 and the fact that 
$D$ is countable that there exists a set $X_0\subset X$
with $\sigma (X\setminus X_0)=0$ such that, for every $x\in X_0$
and every $f\in D,$ 
$$\|T(x)f\|=\|Tf(x)\|=\|Te_0(x)\|\|f\|=h(x)\|f\|$$

Hence, for every $x\in X_0,$ either $h(x)=0$ or the operator
$U_2(x)=T(x)/h(x)$ is an isometry from $D$ to $F.$

The operators $U_2(x)$ can be uniquely extended to isometries on
the whole space $L_1(Y_2).$ In fact, given $a\in L_1(Y_2)$ and
a sequence $a_k\rightarrow a, a_k\in D,$ put
$ U_2(x)(a)=\lim_{k\rightarrow \infty} U_2(x)(a_k).$

Further, for any $a\in L_1(Y_2),$

$$\|a_k - a\|^p = \int_X \|Ta_k(x) - Ta(x)\|^p d\sigma (x) =$$

$$\int_{h(x)=0} \| Ta(x)\|^p d\sigma (x) +
\int_{h(x)\neq 0} \|h(x)U_2(x)a_k - Ta(x)\|^p d\sigma (x)
\rightarrow 0 \tag{4}$$
as $k\rightarrow \infty.$ Therefore, $Ta(x)=0$ for almost
all $x\in X$ with $h(x)=0,$ and $Ta(x)=h(x)U_2(x)a$ for almost
all $x\in X$ (if $h(x)=0$ we put $U_2(x)=0.)$

Similarly, for almost all $x\in X,$ we find isometries $U_1(x)$
from $L_1(Y_1)$ to $F$ such that  $Tb(x)=h(x)U_1(x)b$ for every
$b\in L_1(Y_2).$

Consider an arbitrary function $f\in L_1(Y).$ This function can 
be uniquely represented as a sum $f = f_1 + f_2$ of functions 
$f_1\in L_1(Y_1)$ and $f_2\in L_1(Y_2).$ For all $x\in X$ with
$h(x)\neq 0,$ define operators $U(x)$ 
from $L_1(Y)$ to $F$ by $U(x)f = U_1(x)f_1 + U_2(x)f_2.$ 
By the statement (ii) of Lemma 1, 
for almost every $x$ with $h(x)\neq 0,$ 

$\|U(x)f\|=(1/h(x))\|Tf_1(x) + Tf_2(x)\| = 
\|U_1(x)f_1\| + \|U_2(x)f_2\| = \|f_1\| + \|f_2\| = \|f\|$

Thus, operators $U(x)$ are isometries for almost all $x\in X$
with $h(x)\neq 0.$ In particular, the set $I(L_1,F)$ is
non-empty. Fix any $U\in I(L_1,F)$ and put $U(x)=U$ for
every $x$ with $h(x)=0.$ 

To prove the second statement of Theorem 1 note that, 
for every $f\in L_1(Y),$ we have
$Tf(x) = Tf_1(x) + Tf_2(x) = h(x)( U_1(x)f_1 + U_2(x)f_2)=
h(x)U(x)f$ for almost all $x\in X,$ so $Tf$ and $h(x)U(x)f$
are equal elements of the space $L_p(X,F).$

Now let $E$ be a subspace of $C(K)$ containing the function 
$1(k)\equiv 1$. 

Any function $e\in C(K)$ has a common tangent functional
either with the function $1$ or with the function $-1.$ 
Setting, correspondingly,
$f=1$ or $f=-1$ in Lemma 1 we obtain that, for an arbitrary $e\in E,$
$\|Te(x)\|=\|T1(x)\| \|e\|$ for almost all $x\in X.$
Let $h(x)=\|T1(x)\|.$

Let $e_k, k\in N$ be a sequence of linearly independent functions
from $E$ such that their linear span is dense in $E$ and denote 
by $D$ the set of linear combinations of functions $e_k$ with
rational coefficients. Given fixed representatives $Te_k$ from
the corresponding equivalence classes of functions from 
the space $L_p(X;F),$ define an operator $T(x)$ on $D$ for every 
$x\in X$ by $T(x)(\sum \lambda_{i} e_i)= \sum \lambda_{i} 
Te_{i}(x), \lambda_i \in R.$

Since the set $D$ is countable there exists a set $X_0\subset X$
with $\sigma (X\setminus X_0)=0$ such that, for every $x\in X_0$
and every $e\in D,$ $\|T(x)e\|=\|Te(x)\|=h(x)\|e\|.$ Hence, 
for every $x\in X_0,$ either $h(x)=0$ or the operator
$U(x)=T(x)/h(x)$ is an isometry from $D$ to $F.$

The operators $U(x)$ can be uniquely extended to isometries on
the whole space $E,$ therefore, we have proved the first 
statement of Theorem 1. 

Now an argument similar to (4) proves the second statement.  
\bull \enddemo

\subheading{Remark } If $p=1$ the statement of Theorem 1 is
not true. For instance, the two-dimensional space $l_{\infty}^2$ is
isometric to a subspace of $L_1([0,1]).$ Thus, for any Banach space 
$F,$ $l_{\infty}^2$ is isometric to a subspace of $L_{1}([0,1];F),$  
and Theorem 1 would have implied that $l_{\infty}^2$ is isometric 
to a subspace of any Banach space.

\newpage

\subheading{References}

\item{1.} Bourgain, J.: An averaging result for $c_0$-sequences,
Bull. Soc. Math. Belgique, Vol.B   30(1978), 83-87

\item{2.} Emmanuele, G.: Copies of $l_{\infty}$ in Kothe spaces
of vector valued functions, Illinois J. Math. 36(1992), 293-296 

\item{3.} Koldobsky, A.: Isometries of $L_p(X;L_q)$ and 
equimeasurability, Indiana Univ. Math. J.  40(1991), 677-705

\item{4.} Koldobsky, A.: Measures on spaces of operators and
isometries, J. Soviet Math.  42(1988), 1628-1636

\item{5.} Kwapien, S.: On Banach spaces containing $c_0,$
Studia Math.  52(1974), 187-188

\item{6.} Mendoza, J.: Copies of $l_{\infty}$ in $L^p (\mu ;X),$
Proc. Amer. Math. Soc. 109(1990), 125-127

\item{7.} Pisier, G.: Une propriete de stabilite de la classe
des espaces ne contenant pas $l_1,$ C.R. Acad. Sci. Paris Ser A
86(1978), 747-749

\item{8.} Raynaud, Y.: Sous espaces $l^r$ et geometrie des espaces
$L^p(L^q)$ et $L^{\phi},$ C.R. Acad. Sci. Paris, Ser.I   
301(1985), 299-302

\item{9.} Saab, E. and Saab, P.: On stability problems of some 
properties in Banach spaces, Function spaces (K.Jarosz, editor),
Lect. Notes in Pure and Appl. Math.  136(1992), 367-394

\bye